\documentclass{amsart}
\usepackage{geometry}                % See geometry.pdf to learn the layout options. There are lots.
\usepackage{graphicx}
\usepackage{amssymb}
\usepackage{amsthm}
\usepackage{epstopdf}

\usepackage{url}

\newtheorem{definition}[equation]{Definition}
\newtheorem{defn}[equation]{Definition}

\newtheorem{cor}[equation]{Corollary}
\newtheorem{lem}[equation]{Lemma}
\newtheorem{lemma}[equation]{Lemma}
\newtheorem{addendum}[equation]{Addendum}

\newtheorem{prop}[equation]{Proposition}

\newtheorem{example}[equation]{Example}
\newtheorem{thm}[equation]{Theorem}
\newtheorem{theorem}[equation]{Theorem}
\newtheorem{remark}[equation]{Remark}
\newtheorem{rem}[equation]{Remark}

\newtheorem{question}[equation]{Question}

\def\ol{\overline}
\def\N{\mathbb N}
\def\R{\mathbb R}

\def\Z{\mathbb Z}

\def\Ga{\Gamma}

\def\al{\alpha}
\def\eps{\epsilon}
\def\id{\operatorname{id}}

\def\int{\operatorname{int}}

\def\la{\lambda}
\def\La{\Lambda}

\def\card{\operatorname{card}}

\newenvironment{dedication}
  {      % we want a new page
   \thispagestyle{empty}% no header and footer
   
   \itshape             % the text is in italics
   \raggedleft          % flush to the right margin
  }
  {\par % end the paragraph
  }

\begin{document}

\title{A note on properly discontinuous actions}
\author{M. Kapovich}
\address{Department of Mathematics, UC Davis, One Shields Avenue, Davis CA 95616, USA}
\email{kapovich@math.ucdavis.edu}

\date{\today}                                           

%\subjclass[2020]{37B05, 57S30}

\begin{abstract}
We compare various notions of proper discontinuity for group actions. We also discuss fundamental domains and criteria for cocompactness.  
\end{abstract}

\maketitle

\begin{dedication}
To the memory of Sasha Anan'in
\end{dedication}

\section{Introduction}

This note  is meant to clarify the relation between different commonly used definitions of proper discontinuity without  the local compactness assumption for the underlying topological space. Much of the discussion applies to actions of nondiscrete locally compact Hausdorff topological groups, but, since my primary interest is geometric group theory, I will mostly work  with discrete groups. All group actions are assumed to be continuous, in other words, for discrete groups, these are homomorphisms from abstract groups to groups of homeomorphisms of topological spaces. This combination of {\em continuous} and {\em properly discontinuous}, 
sadly, leads to the ugly terminology ``a continuous properly discontinuous action." 
A better terminology might be  that of a {\em properly discrete} action, since it refers to proper actions of discrete groups. 

Throughout this note, I will be working only with topological spaces which are 1st countable, since spaces most common in metric geometry, geometric topology, algebraic topology and geometric group theory satisfy this property. One advantage of this assumption is that if $(x_n)$ is a sequence converging to a point $x\in X$, then the subset $\{x\}\cup \{x_n: n\in \N\}$ is compact, which is not true if we work with nets instead of sequences. However, I will try to avoid the local compactness assumption whenever possible, since many spaces appearing in metric geometry and geometric group theory (e.g. asymptotic cones) and algebraic topology (e.g. CW complexes) are not locally compact. (Recall that topological space $X$ is {\em locally compact} if every point has a basis of topology consisting of relatively compact subsets.) 

In the last three sections of the note I discuss several concepts related to properly discontinuous actions. 
In Section \ref{sec:cocompact} I discuss cocompactness of group actions. In Section \ref{sec:inmetrics} I discuss group-invariant metrics. In particular, under suitable assumptions I will prove existence of an invariant complete geodesic metric (Theorem \ref{thm:geodesic-metrization}). In Section \ref{sec:fundamentals} I discuss fundamental sets and regions. The main result of this section is Theorem \ref{thm:region} which uses Voronoi tessellations to establish existence of fundamental regions and domains for free properly discontinuous actions on proper geodesic metric spaces.

\medskip
{\bf Acknowledgement.} I am grateful to Boris Okun for pointing out several typos and the reference to \cite{Palais}. I am also grateful to the referee of the paper for useful suggestions and corrections. Lastly, I am grateful to Dennis Schmeckpeper who found a serious mistake in an earlier version of the paper which has lead me to substantially rewrite Section \ref{sec:fundamentals}. 

\medskip 
%{\bf Compliance with Springer Verlag Ethical Standards.} The author has no conflicts of interest to disclose. 

\tableofcontents

\section{Group actions}

A {\em topological group} is a group $G$ equipped with a topology such that the multiplication and inversion maps  
$$
G\times G\to G, (g,h)\mapsto gh, G\to G, g\mapsto g^{-1}
$$
are both continuous. A {\em discrete group} is a group with discrete topology. Every discrete group is clearly a topological group. 

A {\em left continuous action} of a topological group $G$ on a topological space $X$ is a continuous map
$$
\la: G\times X\to X
$$
satisfying 

1. $\la(1_G, x)=x$ for all $x\in X$. 

2. $\la(gh, x)= \la(g, \la(h,x))$, for all $x\in X$, $g, h\in G$.

From this, it follows that the map $\rho: G\to Homeo(X)$ 
$$
\rho(g)(x)= \la(g,x),
$$
is a group homomorphism, where the group operation $\phi \psi$ on $Homeo(X)$ is the composition $\phi\circ \psi$. 
 
If $G$ is discrete, then every homomorphism $G\to Homeo(X)$ defines a left continuous action of $G$ on $X$. 

The shorthand for $\rho(g)(x)$ is $gx$ or $g\cdot x$. Similarly, for a subset $A\subset X$, $GA$ or $G\cdot A$, denotes the {\em orbit  of $A$} under the $G$-action:
$$
GA= \bigcup_{g\in G} gA. 
$$
  
 \medskip 
 The {\em quotient space} $X/G$ (also frequently denoted $G\backslash X$), of $X$ by the $G$-action, is the set of $G$-orbits of points in $X$, equipped with the quotient topology: The elements of $X/G$ are equivalence classes in $X$, where $x\sim y$ when $Gx=Gy$ (equivalently, $y\in Gx$).

The {\em stabilizer} of a point $x\in X$ under the $G$-action is the subgroup $G_x<G$ given by
$$
\{g\in G: gx=x\}. 
$$
An action of $G$ on $X$ is called {\em free} if $G_x=\{1\}$ for all $x\in X$. Assuming that $X$ is Hausdorff, $G_x$ is closed in $G$ for every $x\in X$.

\begin{example}
An example of a left action of $G$ is the action of $G$ on itself via left multiplication:
$$
\la(g, h)= gh. 
$$ 
In this case, the common notation for $\rho(g)$ is $L_g$.  This action is free. 
\end{example}

\section{Proper maps} 

Properness of certain maps is the most common form of defining proper discontinuity; sadly, there are two competing notions of properness in the literature.  

A continuous map $f: X\to Y$ of topological spaces is  {\em proper in the sense of Bourbaki}, or simply {\em Bourbaki--proper} (cf. \cite[Ch. I, \S 10, Theorem 1]{Bourbaki}) if $f$ is a closed map (images of closed subsets are closed) and point--preimages $f^{-1}(y), y\in Y$, are compact. A continuous map $f: X\to Y$ is {\em proper} (and this is the most common definition) if for every compact subset $K\subset X$, $f^{-1}(K)$ is compact. It is noted in \cite[Ch. I, \S 10; Prop. 7]{Bourbaki} that if $X$ is Hausdorff and $Y$ is locally compact, then $f$ is Bourbaki--proper if and only if $f$ is proper. 

The advantage of the notion of Bourbaki-properness is that it applies in the case of Zariski topology, where spaces tend to be compact\footnote{quasicompact in the Bourbaki terminology}  (every subset of a finite-dimensional affine space is Zariski-compact) and, hence, the standard notion of properness is useless. 

Since our goal is to trade local compactness for 1st countability, I will prove a lemma which appears as a corollary in \cite{Palais}: 

\begin{lemma}
If $f: X\to Y$ is proper,  and $X, Y$ are Hausdorff and 1st countable, then $f$ is Bourbaki-proper. 
\end{lemma}
\proof We only have to verify that $f$ is closed. Suppose that $A\subset X$ is a closed subset. Since $Y$ is 1st countable, it suffices to show that for each sequence $(x_n)$ in $A$ such that $(f(x_n))$ 
 converges to $y\in Y$, there is a subsequence $(x_{n_k})$ which converges to some $x\in A$ such that $f(x)=y$. The subset $C=\{y\}\cup \{f(x_n): n\in \N\}\subset Y$ is compact. Hence, by properness of $f$, $K=f^{-1}(C)$ is also compact. Since $X$ is Hausdorff, and $K$ is compact, follows that $(x_n)$ subconverges to a point $x\in K$. By continuity of $f$, $f(x)=y$. Since $A$ is closed, $x\in A$. \qed 

\begin{remark}
This lemma still holds if one were to replace the assumption that $X$ is 1st countable by surjectivity of $f$, see \cite{Palais}. 
\end{remark}

\medskip
The converse (each Bourbaki--proper map is proper) is proven in 
\cite[Ch. I, \S 10; Prop. 6]{Bourbaki} without any restrictions  on $X, Y$. Hence:

\begin{cor}\label{C1}
For maps between 1st countable Hausdorff spaces, Bourbaki-properness is equivalent to properness. 
\end{cor}

\section{Proper discontinuity}

Suppose that $X$ is a 1st countable Hausdorff topological space, $G$ a discrete group and $G\times X\to X$ a (continuous) action. I use the notation  $g_n\to\infty$ in $G$ to indicate that $g_n$ converges to $\infty$ in  the  1-point compactification $G\cup \{\infty\}$ of $G$, i.e. for every finite subset $F\subset G$, 
$$
\card(\{n: g_n\in F\})<\infty. 
$$

Given a group action $G\times X\to X$ and two subsets $A, B\subset X$, the {\em transporter subset} $(A | B)_G$ is defined as
$$
(A | B)_G:=\{g\in G: gA\cap B\ne \emptyset\}. 
$$
Properness of group actions is (typically) stated using certain transporter sets.

\begin{definition}
Two points $x, y\in X$ are said to be {\em $G$-dynamically related} if there is a sequence $g_n\to \infty$ in $G$ and a sequence $x_n\to x$ in $X$ such  that $g_n x_n\to y$. 
\end{definition}

A point $x\in X$ is said to be a {\em wandering point} of the $G$-action if there is a neighborhood $U$ of $x$ such that $(U|U)_G$ is finite. 

\begin{lem}
Suppose that the action $G\times X\to X$ is {wandering} at a point $x\in X$. Then the $G$-action has a {\em $G$-slice} at $x$, i.e. a neighborhood $W_x\subset U$ which is $G_x$-stable and for all $g\notin G_x$, $gW_x\cap W_x=\emptyset$. 
\end{lem}
\proof For each $g\in (U|U)_G - G_x$ we pick a neighborhood $V_g\subset U$ of $x$ such that 
$g V_g\cap V_g= \emptyset$. Then the intersection 
$$
V:=\bigcap_{g\in (U|U)_G - G_x} V_g
$$
satisfies the property that $(V|V)_G= G_x$. Lastly, take 
$$
W_x:= \bigcap_{g\in G_x} V. \qed
$$
%\qed 

The next lemma is clear: 

\begin{lemma}
Assuming that $X$ is Hausdorff and 1st countable, the action $G\times X\to X$ is wandering at $x$ if and only if $x$ is not dynamically related to itself. 
\end{lemma}

Given a group action $\al: G\times X\to X$, we have the natural map 
$$
\hat\al:=\al\times \id_X: G\times X\to X\times X$$ 
where $\id_X: (g,x)\mapsto x$. 

\begin{definition}
An action $\al$ of a discrete group $G$ on a topological space $X$ is {\em Bourbaki--proper} if the map $\hat\al$ is Bourbaki-proper. 
\end{definition}

\begin{lemma}\label{lem:Haus-quotient}
If the action $\al: G\times X\to X$ of a discrete group $G$ on a Hausdorff topological space $X$ is Bourbaki-proper, then the quotient space $X/G$ is Hausdorff. 
\end{lemma}
\proof The quotient map $X\to X/G$ is an open map by the definition of the quotient topology on $X/G$. Since $\al$ is Bourbaki-proper, the image of the map $\hat{\al}$ is closed in $X\times X$. This image is the equivalence relation on $X\times X$ which use used to form the quotient $X/G$. Now,  
Hausdorffness of $X/G$ follows from \cite[Proposition 8 in I.8.3]{Bourbaki}.  \qed 
%Theorem 7.7 in \cite{Tu}.

\begin{definition}
An action $\al$ of a discrete group $G$ on a topological space $X$ is {\em proper} 
if the map $\hat\al$ is proper. 
\end{definition}

Note that the equivalence of (1) and (5) in the following theorem is proven in  \cite[Ch. III, \S 4.4, Proposition 7]{Bourbaki} without any  assumptions on $X$.

\begin{theorem} Assuming that $X$ is Hausdorff and 1st countable, 
the following are equivalent:

\begin{enumerate}
\item[(1)] The action $\al: G\times X\to X$ is Bourbaki-proper.

\item[(2)] For every compact subset $K\subset X$, $$\card( (K|K)_G) <\infty.$$

\item[(3)] The action $\al: G\times X\to X$ is proper, i.e. the map $\hat\al$ is proper.

\item[(4)] For every compact subset $K\subset X$, there exists an open neighborhood $U$ of $K$ such that 
$\card((U|U)_G)<\infty$.  

\item[(5)] For any pair of points $x, y\in X$ there is a pair of neighborhoods $U_x, V_x$ (of $x, y$ respectively) such that $\card((U_x| V_y)_G))<\infty$.  

\item[(6)] There are no $G$-dynamically related points in $X$. 

\item[(7)] Assuming, that $G$ is countable and $X$ is completely metrizable\footnote{It suffices to assume that $X$ is {\em hereditarily Baire}: Every closed subset of $X$ is Baire.}
: The $G$-stabilizer of every $x\in X$ is finite and for any two points $x\in X, y\in X - Gx$, there exists a pair of neighborhoods $U_x, V_y$ (of $x$, resp. $y$) 
such that $\forall g\in G$, $gU_x\cap V_y=\emptyset$. 

\item[(8)] Assuming that $X$ is a metric space and the action $G\times X\to X$ is equicontinuous\footnote{E.g. an isometric action.}: There is no $x\in X$ and a sequence $h_n\to\infty$ in $G$ such that $h_n x\to x$.   

\item[(9)] Assuming that $X$ is a metric space and the action $G\times X\to X$ is equicontinuous: Every $x\in X$ is a wandering point of the $G$-action. 

\item[(10)] Assuming that $X$ is a CW complex and the action $G\times X\to X$ is cellular: Every point of $X$ is wandering. 

\item[(11)] Assuming that $X$ is a CW complex the action $G\times X\to X$ is cellular: Every 
cell in $X$ has finite $G$-stabilizer. 
\end{enumerate}

\end{theorem}
\proof The action $\al$ is Bourbaki-proper if and only if the map $\hat\al$ is proper (see Corollary \ref{C1}) which is equivalent to the statement that for each compact $K\subset X$, the subset $(K | K)_G\times K$ is compact. Hence, (1)$\iff$(2). 

Assume that (3) holds, i.e. $\alpha$ is proper, equivalently, the map $\hat\alpha$ is proper. 
This means that for each compact $K\subset X$, $\hat\al^{-1}(K\times K)= 
\{(g,x)\in G\times K: x\in K, gx\in K\}$ is compact. This subset is closed in $G\times X$ and projects onto $(K|K)_G$ in the first factor and to the subset 
\begin{equation}%\label{eq:union}
\tag{$\star$}
\bigcup_{g\in (K|K)_G} g^{-1}(K) 
\end{equation}
in the second factor. Hence, properness of the action  $\al$ implies finiteness of  $(K|K)_G$, i.e. (2). 
Conversely, if $(K|K)_G$ is finite, compactness of $g^{-1}(K)$ for every $g\in G$
 implies compactness of the union ($\star$). Thus, (2)$\iff$(3).

In order to show that 
(2)$\Rightarrow$(6), suppose that $x, y$ are $G$-dynamically related points: There exists a sequence $g_n\to\infty$ in $G$ and a sequence $x_n\to x$ such that $g_n(x_n)\to y$. The subset
$$
K= \{x,y\}\cup \{x_n, g_n(x_n): n\in {\mathbb N}\}
$$
is compact. However, $y_n\in g_n(K)\cap K$ for every $n$. A contradiction. 

(6)$\Rightarrow$(5): Suppose that the neighborhoods $U_x, V_y$ do not exist. Let $\{U_n\}_{n\in\N}$, $\{V_n\}_{n\in \N}$ be countable bases at $x, y$ respectively. Then for every $n$ there exists $g_n\in G$, such that $g_n(U_n)\cap V_n\ne\emptyset$ for infinitely many $g_n$'s in $G$. After extraction, $g_n\to\infty$ in $G$. 
This yields points $x_n\in U_n, y_n=g_n(x_n)\in V_n$. Hence, $x_n\to x, y_n\to y$. Thus, $x$ is $G$-dynamically related to $y$. A contradiction. 

(5)$\Rightarrow$(4). Consider a compact $K\subset X$. Then for each $x\in K, y\in K$ there exist neighborhoods $U_x, V_y$ such that $(U_x| V_y)_G$ is finite. The product sets $U_x\times V_y, x, y\in K$ constitute an open cover of $K^2$. By compactness of $K^2$, there exist $x_1,..., x_n, y_1,...,y_m\in K$ such that 
$$
K\subset U_{x_1}\cup ... \cup U_{x_n} 
$$
$$
K\subset V_{y_1}\cup ... \cup V_{y_m} 
$$ 
and for each pair $(x_i,y_j)$, 
$$\card(\{g\in G: gU_{x_i}\cap V_{y_j}\ne \emptyset\})<\infty.$$
Setting 
$$
W:= \bigcup_{i=1}^n U_{x_i}, V:= \bigcup_{j=1}^m V_{y_j},$$
we see that 
$$
\card((W|V)_G)<\infty.$$
Taking $U:=V\cap W$ yields the required subset $U$. 

The implication (4)$\Rightarrow$(2) is immediate. 

This concludes the proof of equivalence of the properties (1)---(6). 

\bigskip 
(5)$\Rightarrow$(7): Finiteness of  $G$-stabilizers of points in $X$ is clear. Let $x, y$ be points in distinct $G$-orbits. Let $U'_x, V'_y$ be neighborhoods of $x, y$ such that $(U'_x| V'_y)_G=\{g_1,...,g_n\}$. For each $i$, since $X$ is Hausdorff, there are disjoint neighborhoods $V_i$  of $y$ and $W_i$ of $g_i(x_i)$. Now set
$$
V_y:= \bigcap_{i=1}^n V_i, \quad U_x:= \bigcap_{i=1}^n g_i^{-1}(W_i). 
$$ 
Then $gU_x\cap V_y=\emptyset$ for every $g\in G$. 

(7)$\Rightarrow$(6): It is clear that (7) implies that there are no dynamically related points with distinct $G$-orbits. In particular, every $G$-orbit in $X$ is closed. 

Assume now that $X$ is completely metrizable and $G$ is countable. 
Suppose that a point $x\in X$ is $G$-dynamically related to itself. Since the stabilizer $G_x$ is finite, the point $x$ is an accumulation point of $Gx$; moreover, $Gx$ is closed in $X$. Hence, $Gx$ is a closed perfect subset of $X$. 
Since $X$ admits a complete metric, so does its closed subset $Gx$. Thus, for each $g\in G$, the complement $U_g:=Gx - \{gx\}$ is open and dense in $Gx$. By the Baire Category Theorem, the countable intersection
$$
\bigcap_{g\in G} U_g
$$
is dense in $Gx$. However, this intersection is empty. A contradiction. 

It is clear that (6)$\Rightarrow$(8) (without any extra assumptions). 

(8)$\Rightarrow$(6). Suppose that $X$ is a metric space and the $G$-action is equicontinuous. Equicontinuity implies that for each $z\in X$, a sequence $z_n\to z$ and $g_n\in G$, 
$$
g_n z_n\to gz. 
$$

Suppose that there exist a pair of $G$-dynamically related points $x, y\in X$: $\exists x_n\to x, g_n\in G$, 
$g_n x_n\to y$. By the equicontinuity of the action, $g_n x\to y$. Since $g_n\to\infty$, there exist subsequences $g_{n_i}\to\infty$ and $g_{m_i}\to\infty$ such that  the products 
$h_i:= g_{n_i}^{-1} g_{m_i}$ are all distinct. Then, by the equicontinuity,  
$$
h_i x\to x. 
$$
A contradiction. 

The implications (5)$\Rightarrow$(9)$\Rightarrow$(8) and (5)$\Rightarrow$(10)$\Rightarrow$(11) are clear.

\medskip
Lastly, let us prove the implication (11)$\Rightarrow$(2). We first observe that every CW complex is Hausdorff and 1st countable. Furthermore, every compact $K\subset X$ intersects only finitely many open cells $e_\la$ in $X$. (Otherwise, picking one point from each  nonempty intersection $K\cap e_\la$ we obtain an infinite closed discrete subset of $K$.) Thus, there exists a finite subset $E:=\{e_\la: \la\in \La\}$ of open cells in $X$ such that for every $g\in (K|K)_G$, 
$gE\cap E\ne\emptyset$. Now, finiteness of   $(K|K)_G$ follows from finiteness of cell-stabilizers in $G$. 
\qed

\medskip
Unfortunately, the property that every point of $X$ is a wandering point is frequently taken as the definition of proper discontinuity for $G$-actions, see e.g. \cite{Hatcher, Munkres}.  Items (8) and (10) in the above theorem provide a (weak) justification for this abuse of terminology. I feel that the better name for such actions is {\em wandering actions}.

\begin{example}\label{ex1}
Consider the action of $G=\Z$ on the punctured affine plane $X={\mathbb R}^2 -\{(0,0)\}$, where the generator of $\Z$ acts via $(x,y) \mapsto (2x, \frac{1}{2}y)$. Then for any $p\in  X$, the $G$-orbit $Gp$ has no accumulation points in $X$. However, any two points $p=(x,0), q=(0,y)\in X$ are dynamically related. Thus, the action of $G$ is not proper. 
\end{example}

This example shows that   the quotient space of a wandering action need not be Hausdorff. 

\begin{lemma}
Suppose that $G\times X\to X$ is a wandering action. Then each $G$-orbit is closed and discrete in $X$. In particular, 
the quotient space $X/G$ is T1. 
\end{lemma}
\proof Suppose that $Gx$ accumulates at a point $y$. Then $Gx\cap W_y$ is nonempty, where $W_y$ is a $G$-slice at $y$. It follows that all points of $Gx\cap W_y$ lie in the same $W_y$-orbit, which implies 
that $Gx\cap W_y=\{y\}$. \qed 

\medskip
There are several reasons to consider proper 
actions of discrete (and, more generally, locally compact) groups;  
one reason is that such each proper action of a discrete group yields an {\em orbi-covering map} 
in the case of smooth group actions on manifolds: $M\to M/G$ is an orbi-covering provided that the action of $G$ on $M$ is smooth (or, at least, locally smoothable). Another reason is that for a proper action on a Hausdorff space, $G\times X\to X$, the quotient $X/G$ is again Hausdorff, see Lemma \ref{lem:Haus-quotient}.  

\begin{question}
Suppose that $G$ is a discrete group, $G\times X\to X$ is a free continuous action on an $n$-dimensional topological manifold $X$ 
such that the quotient space $X/G$ is a (Hausdorff)  $n$-dimensional topological manifold. Does it follow that the $G$-action  on $X$ is proper? 
\end{question}
 
The answer to this question is negative if one merely assumes that $X$ is a locally compact Hausdorff topological space and $X/G$ is Hausdorff, see \cite{MO4} (the action given there was even cocompact). 
Below is a different example.  We begin by constructing a non-proper 
free continuous $\R$-action on a manifold, such that the quotient space is not just Hausdorff but  is a manifold with boundary. 

\begin{example}\label{ex2}
This is a variation on Example \ref{ex1}. We start with the space 
$$
Z=\{(x,y): x, y\in [0,\infty), (x,y)\ne (0,0)\}.
$$
 Take the quotient space $X$ of $Z$ by the equivalence relation $(x,0)\sim (0, \frac{1}{x})$. The space $X$ is homeomorphic to the open Moebius band. The group $G=\R$ acts on $Z$ continuously by
$$
(t, (x,y))\mapsto (2^t x, 2^{-t}y).
$$
The above equivalence relation on $X$ is preserved by the $G$-action and, hence, the $G$-action descends to a continuous $G$-action on $X$. It is easy to see that this action is free but not proper: The equivalence class of $(1,0)$ is dynamically related to itself. Lastly, the quotient $X/G$ is Hausdorff, homeomorphic to 
$[0, 1)$ (the equivalence class of $(1,0)$ maps to $0\in [0,1)$). 
 \end{example}

Lastly, we use Example \ref{ex2} to construct a non-proper free $\Z$-action with Hausdorff quotient. We continue with the notation of the previous example. 

\begin{example}
Let $Y\subset Z$ denote the following subset of $Z$ (with the subspace topology): 
$$
Y=\{(2^m, 0): m\in \Z\} \cup \{(0, 2^n): n\in \Z\} \cup \{(2^m, 2^n): (m,n)\in \Z^2\}. 
$$
Let $W$ denote the projection of $Y$ to $X$. We take $\Gamma=\Z< G=\R$. This subgroup preserves $Y$ and, hence, $W$.  The quotient $W/\Gamma$ is homeomorphic to $Y\cap \{(0, y): y\in \R\}$, hence, is Hausdorff.   At the same time, the $\Gamma$-action on $W$ is  non-proper. 
\end{example}

\section{Cocompactness}\label{sec:cocompact}

There are two common notions of cocompactness for group actions: 

\begin{enumerate}

\item $G\times X\to X$ is cocompact if there exists a compact $K\subset X$ such that $G\cdot K=X$.

\item $G\times X\to X$ is cocompact if $X/G$ is compact. 

\end{enumerate}

It is clear that (1)$\Rightarrow$(2), as the image of a compact under the continuous (quotient) map $p: X\to X/G$ 
is compact. 

\begin{lemma}
If $X$ is locally compact then (2)$\Rightarrow$(1). 
\end{lemma}
\proof For each $x\in X$ let $U_x$ denote a relatively compact neighborhood of $x$ in $X$. Then 
$$
V_x:=p(U_x)= p(G\cdot U_x),
$$
is compact since $G\cdot U_x$ is open in $X$. Thus, we obtain an open cover $\{V_x: x\in X\}$ of $X/G$. Since $X/G$ is compact, this open cover contains a finite subcover 
$$
V_{x_1},..., V_{x_n}.
$$
It follows that 
$$
p( \bigcup_{i=1}^n {U_{x_i}} )=X/G. 
$$
The set 
$$
K= \bigcup_{i=1}^n \ol{U_{x_i}}
$$
is compact and $p(K)=X/G$. Hence, $G\cdot K=X$. \qed 

\begin{lemma}
Suppose that $X$ is normal and Hausdorff, $G\times X\to X$ is a proper action of a discrete group, such that $X/G$ is locally compact. Then $X$ is locally compact. 
\end{lemma}
\proof Pick $x\in X$. Let $W_x$ be a slice for the $G$-action at $x$; then $W_x/G_x\to X/G$ is a topological embedding. Thus, our assumptions imply that $W_x/G_x$ is compact for every $x\in X$. 
Let $(x_\al)$ be a net in $W_x$. Since $W_x/G_x$ is compact, the net $(x_\al)/G$ contains a convergent subnet. Thus, after passing to a subnet, there exists $g\in G_x$ such that  $(gx_\al)$ converges to some 
$x\in \ol{W_x}$.  Hence, $(x_\al)$ subconverges to $g^{-1}(x)$. Thus, $W_x$ is relatively compact. Since $X$ is assumed to be normal, $x$ admits a basis of relatively compact neighborhoods. \qed  

\begin{cor}
For normal Hausdorff spaces $X$ the two notions of cocompactness agree for proper discrete group actions on $X$.  
\end{cor}

On the other hand, if the drop the properness condition, the two notions are not equivalent even for $\Z$-actions with Hausdorff quotients, see the example by R. de la Vega in \cite{MO5}.

\section{Invariant metrics}\label{sec:inmetrics}

We start with several general definitions. 
A discrete subset $E$ of a metric space $(X,d)$ will be called {\em metrically proper} if for some (equivalently, every) $p\in X$ the function
$$
d(p, \cdot): E\to \R_+
$$
is proper. In other words, every metric ball contains only finitely many points of $E$. A {\em geodesic metric space}, is a metric space $(X,d)$ where every two points $x, y$ are connected by a geodesic segment, i.e. an isometric embedding $c: [a,b]\to (X,d)$ such that $c(a)=x, c(b)=y$. Geodesic segments connecting $x$ to $y$ need not be unique; however, one frequently denotes such segments $xy$ by 
abusing the notation. We will also conflate geodesic segments and their images. Note that each locally compact complete geodesic metric space $(X,d)$ is {\em proper}, i.e. closed metric balls in $(X,d)$ are compact, 
see \cite[Theorem 2.5.28]{BBI}.

An isometric action $G\times X\to X$ of a discrete group is {\em metrically proper} if $G$ acts with finite point-stabilizers and one (equivalently, every) $G$-orbit in $X$ is a metrically proper subset. In other words, 
for every $x\in X$ the function 
$$
g\mapsto d(x, gx)
$$
is proper on $G$. This condition is stronger than properness of the action but is equivalent to properness of the $G$-action in the case of proper metric spaces $(X,d)$. Given an isometric properly discontinuous $G$-action on $X$ we define the function
$$
\rho: X/G\to \R_+
$$
sending each equivalence class $[x]\in X/G$ to 
$$
\inf\{d(gx,x): g\in G\setminus \{1\}\}. 
$$
This function is $2$-Lipschitz:  
$$
|\rho([x]) - \rho([y])|\le 2 d([x], [y]). 
$$
If the $G$-action is metrically proper, then the infimum in the definition of $\rho$ is realized and if the action is also free then $\rho([x])>0$ for all $x\in X$. By abusing the notation, we will also denote this function $\rho(x)$. 

\medskip 
Suppose that $(X,d)$ is a metric space and $G$ is a group acting isometrically  and metrically properly on $X$. One defines the {\em quotient-metric} $d_G$ on $X/\Gamma$ by
\begin{equation}\label{eq:q-metric} 
d_G([x], [y])=\min_{g\in G} d(x, G y)= \min_{g, h\in G} d(gx, hy),
\end{equation}
where $[x], [y]\in X/G$ are equivalence classes of points $x, y\in X$ under the equivalence relation defined by $G$. Then $d_G$ is a metric on $X/G$ which metrizes the quotient topology on $X/G$, see \cite[Theorem 6.6.2]{Ratcliffe}. By the construction, the quotient map $q: (X,d)\to (X/G, d_G)$ is $1$-Lipschitz. 

\begin{lem}\label{lem:quot-space}
Suppose that the $G$-action on $X$ is metrically proper. Then the following hold:

1. If the metric space $(X,d)$ is geodesic and complete, then so is $(X/G, d_G)$.  

2. If $(X,d)$ is proper, so is $(X/G, d_G)$. 

3. If the $G$-action is free then the quotient map $q: (X,d)\to (X/G, d_G)$ is a local isometry. More precisely, for every $x\in X$ the restriction of $q$ to $B(x, \frac{1}{8}\rho(x))$ is an isometry onto $B([x], \frac{1}{8}\rho(x))$. 
\end{lem}
\proof 1a. Take points $[x], [y]\in X/G$. Pick their representatives $x, y\in X$ which realize the minimal distance between the corresponding $G$-orbits in $X$. Let $c: [0,T]\to xy\subset X$ be a geodesic connecting $x$ to $y$. Then, by the definition of the metric $d_G$, the composition of $c$ with the quotient map $q: X\to X/G$ is a geodesic in  $(X/G, d_G)$ connecting  $[x]$ to $[y]$.

1b. Suppose that $(z_n)$ is a Cauchy sequence in $X/G$. Then the diameter $D$ of the subset 
$$\{z_n: n\in \N\}\subset X/G$$
 is finite. We inductively choose a subsequence $(z_{n_i})$ in $(z_n)$ such that 
$$
d_G(z_{n_i}, z_{n_{i+1}})\le\frac{D}{ 2^i}, i\in \N. 
$$
Concatenating geodesic segments $z_{n_i} z_{n_{i+1}}$  we obtain a piecewise-geodesic path 
$$
\gamma: [0, T)\to X
$$ 
whose length $T$ is at most 
$$
\sum_{n=1}^\infty \frac{D}{ 2^i}<\infty. 
$$
We then inductively lift each geodesic segment in $\gamma$ to a geodesic segment in $X$ and obtain 
a piecewise-geodesic path $c: [0, T)\to X$ of length $T$. Since $(X,d)$ is complete, the path $c$ extends continuously to $T$. 
Projecting $c(T)$ to $X/G$ we obtain the limit of the subsequence $(z_{n_i})$. Hence, $(z_n)$ converges as well. 

2. Suppose that $(X,d)$ is proper. Consider the closed metric ball $\bar{B}([x], R)$ in $(X/G, d_G)$. Then the closed ball 
$\bar{B}(x, R)\subset X$ projects onto  $\bar{B}([x], R)$. Compactness of  $\bar{B}(x, R)$ implies compactness of  $\bar{B}([x], R)$. 

3. Fix $x\in X$, set $R=\frac{1}{8}\rho(x)$ and consider points $y, z\in B(x, R)$. We have to verify that $d(y,z)=d(y, Gz)$. Take 
$g\in G\setminus \{1\}$. We have $|\rho(y)-\rho(x)|\le 2R$ and $|\rho(z)-\rho(x)|< 2R$ since $\rho$ is $2$-Lipschitz. 
Thus, $d(z, gz)> \rho(x)- 2R= \rho(x) - \frac{1}{4}\rho(x)= \frac{3}{4}\rho(x)$. By the triangle inequality,
$$
d(y, gz) >  \frac{3}{4}\rho[x)- 2R=  \frac{3}{4}\rho(x) - \frac{1}{4}\rho(x)=  \frac{1}{2}\rho(x) > 2R> d(y,z). 
$$
Lastly, for every $r>0$ and $x\in X$, $q(B(x,r))$ is contained in $B([x],r)$ since the quotient map $q: (X,d)\to (X/G, d_G)$ is $1$-Lipschitz. 
In the case $r= R$ as above, the fact that $q$ restricts to an isometry on $B(x,R)$ implies the equality $q( B(x, R))= B([x], R)$. 
\qed

\medskip 
It turns out that under some rather mild assumptions, given a proper action  $G\times X\to X$, 
there is a $G$-invariant metric metrizing the topology on $X$:

\begin{thm}\label{thm:inmetric}
Suppose that $G$ is a locally compact Hausdorff group,  $X$ is locally compact, metrizable space, $G\times X\to X$ is a proper action  and $X/G$ is paracompact. Then $X$ admits a $G$-invariant metric metrizing the topology on $X$
\end{thm}

See \cite[Theorem 3]{Koszul}. Koszul also notes that if $X$ is paracompact and locally connected, then $X/G$ is paracompact. This theorem was improved in \cite{AMN}:

\begin{thm}
Suppose that $G$ is a locally compact Hausdorff group,  $X$ is locally compact, $\sigma$-compact metrizable space, and $G\times X\to X$ is a proper action. Then $X$ admits a $G$-invariant proper metric metrizing the topology on $X$. 
\end{thm}

A  Riemannian version of these theorems holds in the context of smooth actions of Lie groups: 

\begin{thm}
Suppose that $X$ is a smooth manifold, $G$ is a Lie group and $G\times X\to X$ is a smooth proper action. There there exists a $G$-invariant complete Riemannian metric on $X$. 
\end{thm}

See \cite[Theorem 2]{Koszul} for the existence of an invariant Riemannian metric  and  \cite{Kankaanrinta}  for the existence of an invariant complete Riemannian metric. 

\medskip 
We next discuss a construction of $G$-invariant {\em complete geodesic metrics} on more general topological spaces.

\begin{theorem}\label{thm:geodesic-metrization} 
Suppose that $X$ is a 2nd countable, connected and locally connected locally compact Hausdorff topological space. 
Suppose that $G\times X\to X$ is a proper action of a discrete countable group such that the fixed-point set of each nontrivial element of $G$ is nowhere dense in $X$. Then $X$ can be metrized using a $G$-invariant complete geodesic metric. 
\end{theorem}  
 \proof  
 
 \begin{lemma}%\label{lem:quot-space} 
 The quotient space $Y=X/G$ is locally compact, connected, locally connected and metrizable.
 \end{lemma} 
 \proof Local compactness and connectedness of $Y$ follows from that of $X$. The 2nd countability of $X$ implies the 2nd countability of $Y$. By Lemma \ref{lem:Haus-quotient}, $Y$ is Hausdorff. Since $Y$ is locally compact and Hausdorff, its one-point compactification is compact and Hausdorff, hence, regular. It follows that $Y$ itself is regular. In view of the 2nd countability of $Y$, Urysohn's metrization theorem implies that $Y$ is metrizable. \qed
 
 \begin{rem}
 Note that each locally compact metrizable space is also locally path-connected.
 \end{rem}

It is proven in \cite{TT} that each locally compact, connected, locally connected metrizable space, such as $Y$,  admits a complete geodesic  metric $d_Y$ which we fix from now on. 
Consider the projection $p: X\to Y$.  According to \cite[Theorem 6.2]{Bredon} (see also \cite[Lemma 2]{A}), the map $p$ satisfies the path-lifting property: Given any path $c: [0,1]\to Y$, a point $x\in X$ satisfying $p(x)=c(0)$, there exists a path $\tilde{c}: [0,1]\to X$ such that $p\circ \tilde{c}=c$. (This result is, of course, much easier if the $G$-action is free, i.e. $p: X\to Y$ is a covering map.) 
We let ${\mathcal L}_X$ denote the set of paths in $X$ which are lifts of rectifiable paths $c: [0,1]\to Y$. Clearly, the postcomposition of $\tilde{c}\in {\mathcal L}_X$ with an element of $G$ is again in ${\mathcal L}_X$.  Our next goal is to equip $X$ with   a $G$-invariant {\em length structure} using the family of paths ${\mathcal L}_X$. Such a structure is a function on ${\mathcal L}_X$ with values in $[0,\infty)$, satisfying certain axioms that can be found in \cite[Section 2.1]{BBI}. Verification of most of these axioms is straightforward, I will check only some (items 1, 2, 3 and 4 below). 

1. If $\tilde{c}\in {\mathcal L}_X$ is a lift of a path $c$ in $Y$, then we declare $\ell(\tilde{c})$ to be equal to the length of $c$. 

2. If $\tilde{c}_i, i=1, 2$, are paths in  ${\mathcal L}_X$ (which are lifts of the paths $c_1, c_2$ respectively) whose concatenation $b=\tilde c_1\star \tilde c_2$ is defined, then $b$ is a lift of the concatenation $c_1\star c_2$. Clearly, $\ell(b)= \ell(\tilde{c}_1) +\ell(\tilde{c}_2)$. 

3. Let $U$ be a neighborhood of some $x\in X$. We need to prove that 
\begin{equation}\label{eq:positive}
\inf_{\gamma} \{ \ell(\gamma)\} >0,
\end{equation}
where the infimum is taken over all $\gamma=\tilde{c}\in {\mathcal L}_X$ connecting $x$ to points of $X\setminus U$. It suffices to prove this claim in the case when $U$ is $G_x$-invariant,  satisfies 
\begin{equation}\label{eq:stablenbd}
\ol{U}\cap g\ol{U}\ne \emptyset \iff g\in G_x, 
\end{equation}
and $\gamma$ connects $x$ to points of $\partial U$. Then $V=p(U)$ is a neighborhood of $y=p(x)$ in $Y$ and the paths $c=p\circ \gamma$ connect $y$ to points in $\partial V$. But the lengths of the paths $c$ are clearly bounded away from zero and are equal to the lengths of their lifts $\tilde c$. Thus, we obtain the required bound \eqref{eq:positive}. 

4. Let us verify that any two points in $X$ are connected by a path in ${\mathcal L}_X$. Since $X$ is connected, it suffices to verify the claim locally. Let $U$ is $G_x$-invariant neighborhood of $x$ satisfying \eqref{eq:stablenbd}, such that $V=p(U)$ is an open metric ball in $Y$ centered at $y=p(x)$. Take $u\in U$, $v:= p(u)\in V$. 
Let $c: [0,T]\to V$ be a geodesic connecting $v$ to $y$. Then there exists a lift $\tilde{c}: [0,T]\to U$ of $c$ with $\tilde{c}(0)=u$. Since $x\in U$ is the only point projecting to $y$, we get 
$\tilde{c}(T)=x$. By taking concatenations of pairs of such radial paths in $U$, we conclude that any two points in $U$ are connected by a path $\tilde{c}\in {\mathcal L}_X$.

\medskip 
Given a length structure on $X$, one defines a path-metric (metrizing the topology of $X$) by
$$
d_X(x_1, x_2)= \inf_\gamma \{\ell(\gamma)\}
$$
where the infimum is taken over all $\gamma\in {\mathcal L}_X$ connecting $x_1$ to $x_2$. By the construction, the projection $p: (X,d_X)\to (Y,d_Y)$ is 1-Lipschitz. 

\begin{lemma}
The metric $d_X$ is complete. 
\end{lemma}
\proof Let $(x_n)$ be a Cauchy sequence in $(X,d_X)$. By the construction of the metric $d_X$, there exists a finite length  
path $\tilde{c}: [0,1)\to (X,d_X)$ and a sequence $t_n\in [0,1)$ such that $\tilde{c}(t_n)=x_n, \tilde{c}(0)=x=x_1$. Since the map $p$ is 1-Lipschitz, the path $c=p\circ \tilde{c}: [0,1)\to (Y,d_Y)$ also has finite length. Since the metric $d_Y$ was complete to begin with, the path $c$ extends to a path $\bar{c}: [0,1]\to Y$; set $y':=\bar{c}(1)$. 

Assume for a moment that $G$ acts freely on $X$. Then we have  the {\em uniqueness} of lifts of paths from 
$Y$ to $X$. Thus, the unique lift $\tilde{\bar{c}}$ of $\bar{c}$ starting at the point $x$ satisfies the property that its restriction to $[0,1)$ equals $\tilde{c}$. It follows that the sequence $(x_n)$ converges to $\tilde{\bar{c}}(1)$. Below we  generalize this argument to the case of non-free actions.

Let $U$ be a neighborhood of $y'=\bar{c}(1)$ which is the projection to $Y$ of a relatively compact 
slice neighborhood $\tilde{U}$ of some 
$x'\in p^{-1}(y')$. Without loss of generality (by removing finitely many initial terms of the sequence $(x_n)$) we can assume that the image of the path $c$ lies entirely in $U$. Applying the path-lifting property to the path $c$ with the prescribed {\em terminal} point $x'$, we obtain a lift of the path  $\bar{c}$ that terminates at $x'$. This lift has to be entirely contained in $\tilde{U}$ and its initial point has to be of the form $g(x)$ for some $g\in G$. Applying $g^{-1}$ to this lift, we obtain another lift of $\bar{c}$, denoted  $\tilde{\bar{c}}$, which starts at $x$ and terminates at $g^{-1}(x')$. 

Consider the restriction of $\tilde{\bar{c}}$ to $[0,1)$. This restriction is also a lift to the path $c|_{[0,1)}$ and the image of the latter lies entirely in $U$. Hence, the image of $\tilde{\bar{c}}|_{[0,1)}$ lies entirely in the relatively compact subset 
$g^{-1}(\tilde{U})\subset X$. Thus, the Cauchy sequence 
$(x_n)$ lies in a relatively compact subset of $X$, and it follows that this sequence converges in $X$. \qed

\medskip 
Since $(X,d_X)$ is locally compact and complete, by Theorem 2.5.28 (and Remark 2.5.29) in \cite{BBI}, $(X,d_X)$ is a geodesic metric space. Lastly, we note that, by the construction, the length structure on $X$ and, hence, the metric $d_X$, is $G$-invariant. This concludes the proof of the theorem. \qed 

\begin{question}
Local compactness and local connectivity were critical for the proof of the theorem. Does the theorem hold without these  assumptions? 
\end{question}

\section{Fundamental domains of properly discontinuous group actions}\label{sec:fundamentals}
  
\subsection{Fundamental sets}\label{sec:fundamental sets}

As with many notions going back to the 19th century, there is no consistency in the literature regarding the definition of fundamental sets and domains.   The next definition follows  \cite{Koszul}. Our definition is similar to the definition given by Borel and Ji in 
\cite[Definition III.2.14]{BJ}, except that their local finiteness condition is weaker: It is required only for singletons $K$.

 \begin{definition}\label{def:fundamental set}
 A closed subset $F\subset X$ is a {\em fundamental set} for a proper action of a discrete $G$ on a topological space 
 $X$ if $G\cdot F=X$ and 
 for every compact $K\subset X$, the transporter set $(F|K)_G$ is finite (the {\em local finiteness condition}). A closed 
subset $F\subset X$ is a {\em fundamental set in the sense of Koszul} if, moreover, 
there exists an open neighborhood $U$ of $F$ such that for every compact $K\subset X$, the transporter set $(U|K)_G$ is finite.  \end{definition} 
  
 Fundamental sets appear naturally in the reduction theory of arithmetic groups (Siegel sets), see \cite{Siegel} and \cite{BJ}. We note, however, that in the literature there are many alternative notions of fundamental sets, inconsistent with the one given above, see e.g. Beardon's book \cite[9.1]{Beardon}: According to Beardon's definition, a subset $F$ of $X$ is called fundamental for the action of $G$ on $X$ if $F$ intersects every $G$-orbit in $X$ in exactly one point.  We will avoid using this definition since its set-theoretic nature provides us with no useful control of the structure of $F$. 
 
 \medskip
The local finiteness condition in the definition of a fundamental set has  several implications:

\begin{lemma}\label{lem:nbd}
Suppose that $F\subset X$ is a fundamental set for a proper action of a discrete group $G$ on a 
1st countable and Hausdorff space $X$. Then:

1.  For every $x\in X$ there exists a neighborhood $W$ of $x$ such that $(F|U)_G$ is finite.  

2. For every $x\in X$ there exist a finite subset $E=\{g_1,...,g_k\}\subset G$ such that the interior of $g_1(F)\cup...\cup g_n(F)$ is a neighborhood of $x$ in $X$. 
\end{lemma}
\proof 1. Suppose that such $W$ does not exist. Then there exists a sequence of distinct elements $g_n\in G$ and points $x_n\in X$ such that 
$$
\lim_{n\to\infty} x_n=x
$$
and $x_n\in g_n(F)$. It follows that for the compact $K=\{x_n: n\in \mathbb N\}\cup \{x\}$ the transporter set $(F|K)_G$ is infinite, which is a contradiction. 

2. By the local finiteness condition, there are only finitely many elements $g_1,...,g_k\in G$ such that $x\in g_i(F)$. By Part 1 of the lemma, there exists a neighborhood $W$ of $x$ such that $W\cap gF\ne \emptyset$ only for $g\in E=\{g_1,...,g_k\}$. But then, since $GF=X$, it follows that
$$
W\subset g_1F\cup ... \cup g_kF. \qed
$$
%\qed

  \medskip 
 For each fundamental set $F$ of a $G$-action on a topological space $X$ we define its quotient space $F/G$ as the quotient space of the equivalence relation 
 $x\sim y$ $\iff$ $Gx=Gy$. %$(\{x\}|\{y\})_G\ne\emptyset$. 
 The following proposition explains why fundamental sets are useful: They allow one to describe quotient spaces of proper actions by discrete groups using less information than is contained in the description of the action.  
  
 \begin{prop}
 Suppose that $F$ is a fundamental set for a proper action by discrete group $G$ 
 on a 1st countable and Hausdorff space $X$. Then the natural projection map $p: F/G\to X/G$ is a homeomorphism. 
 \end{prop} 
 \proof The  map $p$ is continuous by the definition of the quotient topology. It is also obviously a bijection. It remains to show that $p$ is a closed map. Since $F$ is closed, it suffices to show that the projection $q: F\to X/G$ is a closed map. 
 Suppose that $(x_n)$ is a sequence in $F$ 
 such that $q(x_n)$ converges to some $y\in X/G$, $y$ is represented by a point $x\in F$. Then there is a sequence $h_n\in G$ such that $z_n=h_n(x_n)$ converges to $x$. If the sequence $(h_n)$ contains infinitely many distinct elements, we obtain a contradiction with the local finiteness property of $F$ similarly to the proof of Lemma \ref{lem:nbd}. Hence, the set  $E=\{h_n: n\in \mathbb N\}$ is finite. 
 Applying inverses of the elements  
 $h\in E$, to the sequence $(z_n)$, we see that the subset $\{x_n: n\in \mathbb N\}\subset X$ is relatively compact.  
 Thus, $q: F\to F/G$ is a closed map. \qed

 \medskip 
  There are several existence theorems for fundamental sets. The next proposition, proven in \cite[Lemma 2]{Koszul}, guarantees  existence of fundamental sets under the {\em paracompactness assumption} on $X/G$. 
  
\begin{prop}
 Each proper action 
$G\times X\to X$ of a discrete group $G$ on a locally compact Hausdorff space $X$ with paracompact quotient $X/G$ admits  a fundamental set in the sense of Koszul. 
\end{prop}  
  
Another construction of fundamental sets is given by {\em closed Dirichlet domains}. Let $G\times X\to X$  be an isometric  proper action 
of a discrete group $G$ on a metric space $(X,d)$. The {\em closed Dirichlet domain} for this action is 
\begin{equation}\label{eq:hat}
\hat{D}_x=\{y\in X: d(y,x)\le d(y, gx)\quad \forall g\in G\}. 
\end{equation}
Note that $g \hat D_x= \hat{D}_{gx}$. We also note that $\hat{D}_x$ is a closed subset of $X$ since it is the intersection of a family of closed subsets 
$$
\{y\in X: d(y,x)\le d(y, gx)\}, g\in G. 
$$ 

\begin{prop}\label{prop:Dirichlet domain}
Suppose that $G\times X\to X$ is a metrically proper isometric action of a discrete group $G$. Then every closed Dirichlet domain $\hat D= \hat D_x$ is a fundamental set for the $G$-action.
\end{prop}
\proof  1. Let us prove that $g \hat{D}=X$. For each $y\in X$ the function $g\mapsto d(y, gx)$ is a proper function on $G$, hence, it attains its minimum at some $g\in G$. Then, clearly, $y\in \hat{D}_{gx}= g\hat{D}_x$. Thus, $g \hat{D}=X$.  

2. Secondly, we verify local finiteness. Consider a metric ball $B=B(x,R)$ for any $R>0$. If $\hat{D}_{gx}\cap B\ne \emptyset$, for every point $y$ in this intersection 
$$
d(g^{-1}y, x)=d(y, gx)\le d(y, x)< R,
$$ 
In view of metric properness of the $G$-action, the set of such elements $g\in G$ is finite. \qed 

\medskip
We will discuss Dirichlet domains (and their generalizations via Voronoi tessellations) again in Section \ref{sec:fundamental domains}.

\medskip
Note that in the definition of a closed Dirichlet domain one does not really need a metric, what is needed is a $G$-invariant continuous function $d: X\times X\to \mathbb R_+$. For the proof of   Proposition \ref{prop:Dirichlet domain} to go through one needs a metric $\delta$ on $X$ such that:

(a) The $G$-action is metrically proper on $(X,\delta)$. 

(b)   $\delta(y, x)\le \phi(d(y,x))$ for some function $\phi$.

\medskip 
An example of the situation when this is useful appears in the context of discrete subgroups $\Gamma$ of $G=SL(n,\R)$ acting on the space $X$ of symmetric  positive-definite $n\times n$ matrices $M$ with $\det M=1$ by 
$$
M\mapsto g^T M g, g\in SL(n,\R).
$$
Then Selberg in \cite{Selberg} used the function $d: X\times X\to \R_+$,  
$$
d(A,B)= \log\left(\frac{1}{n}tr(A^{-1}B)\right)
$$
to define an analogue of Dirichlet domains for the $\Gamma$-action on $X$. (See also \cite{Kapovich-2023}.) The advantage of such generalized   Dirichlet domains is that they are intersections of $X$ with polyhedral cones in the space of all symmetric $n\times n$ matrices. 
  
 \begin{defn}
 Suppose that $G\times X\to X$ is a continuous action. A closed subset $F\subset X$ is a {\em strict fundamental set} for the action if it intersects each $G$-orbit in $X$ in exactly one point. 
 \end{defn} 
  
Strict fundamental sets do not exist often, but they do exist for some classes simplicial group actions on simplicial complexes (one does not even need to assume properness), e.g. for actions of Coxeter groups on {\em Coxeter complexes} and actions of semisimple Lie groups (as well as semisimple algebraic groups over discrete valued fields) on {\em buildings} (see e.g. \cite{Ronan}).  In the next section we will use a construction of strict fundamental sets for properly discontinuous 
simplicial group actions on vertex sets of connected graphs described below.   
%Below we construct {\em connected} strict fundamental sets for groups acting on connected graphs:
Suppose that $\Gamma$ is a simplicial graph (a 1-dimensional simplicial complex), $G\times \Gamma\to \Gamma$ is a simplicial action of a discrete group $G$. (The action need not be proper.) 
Note that the edge-stabilizers need not fix the invariant edges. However, if $\Gamma'$ denotes the barycentric subdivision of $\Gamma$ then the induced action of $G$ on $\Gamma'$ is {\em without inversions}, i.e. if an element of $G$ preserves an edge, then it fixes the edge pointwise.

\begin{lem}\label{lem:tree}
Suppose that $\Gamma$ is connected. Then there exists a subtree $\Phi\subset \Ga'$ such that the vertex set of $\Phi$ is 
a strict fundamental set  for the $G$-action on the vertex set of $\Ga'$. 
\end{lem}
\proof The quotient $\Ga'/G$ has natural structure of a connected simplicial graph. Let $q: \Ga'\to \Ga'/G$ denote the quotient map. 
Choose $T\subset \Ga'/G$, a maximal subtree (this may require the Axiom of Choice if the vertex set of $\Ga'/G$ is uncountable). We will construct $\Phi$ by lifting $T$ (inductively) to $\Gamma'$. We pick a vertex $v\in \Ga'/G$ and lift it arbitrarily to a vertex $\tilde v\in q^{-1}(v)\in \Ga'$. Then, of course, $G \{\tilde v\}\cap \{\tilde v\}= \{\tilde v\}$. We proceed inductively, working with subtrees $B_n\subset T$ which are closed metric balls of radius $n$ centered at $v$.   
Suppose that we defined a subtree $\Phi_n\subset \Ga'$  such that $q(\Phi_n)=B_n$ and each $G$-orbit in $\Ga'$ intersects $\Phi_n$ in at most one point. Let $e=[u,w]$ be an edge in $B_{n+1}$ with $u\in B_n$. Then there exists an edge $\tilde{e}=[\tilde{u}, \tilde{w}]$ of $\Ga'$ which projects to $e$ and $\tilde{u}\in B_n$ is a vertex projecting to $u$. We add the edge $\tilde{e}$ (and the vertex $\tilde w$) 
to $\Phi_n$ (note that $\tilde w$ cannot belong to $\Phi_n$). We repeat this for all edges of $B_{n+1}$ which are not in $B_n$, resulting in a subtree $\Phi_{n+1}\subset \Ga'$. By the construction, each $G$-orbit in $\Ga'$ intersects  
$\Phi_{n+1}$ in at most one point. Lastly, the union
$$
\Phi=\bigcup_{n} \Phi_n
$$
is a subtree satisfying the required properties. \qed 

\medskip 
Note that, unless $\Ga'/G=T$ (i.e. $\Ga'/G$ is a tree), $\Phi$ is not a fundamental set of the $G$-action on $\Ga'$ since preimages of edges of  $\Ga'/G$ that are not in $T$ are not contained in the $G$-orbit of $\Phi$.

\subsection{Fundamental regions and domains}\label{sec:fundamental domains}
   
One frequently encounters a   sharper version of fundamental sets, called {\em fundamental domains} or {\em fundamental regions}. Again, there is no consistency in this definition in the literature. Below is a small sample of existing definitions. Ratcliffe in \cite[\S 6.6]{Ratcliffe} defines {\em fundamental regions} for a properly discontinuous isometric $G$-action  on a metric space $(X,d)$  as {\em open} subsets $R\subset X$ such that $X=G \ol{R}=X$ and $gR\cap R=\emptyset$ for all $g\in G\setminus \{1\}$. Then Ratcliffe defines {\em fundamental domains} as connected fundamental regions. Ratcliffe also defines {\em locally finite} fundamental domains by imposing the extra assumption of local finiteness just as in Definition \ref{def:fundamental set} given above. Beardon in \cite[\S 9.1, 9.2]{Beardon} also defines fundamental domains as open connected subsets as above, but  (working in the context of subsets of hyperbolic spaces) imposes the extra condition that the boundary has Lebesgue measure zero. In contrast, S.~Katok in \S 3.1 of \cite{Katok}, defines fundamental regions $F\subset X$ as closures of certain open subsets $R\subset X$ where $R$ is a fundamental region as in Ratcliffe's definition. Furthermore, Benedetti and Petronio, \cite[\S C1]{BP}, define fundamental domains as {\em Borel subsets} $F\subset X$ such that $GF=X$ and $gF\cap F\subset \partial F$ for all $g\in G\setminus \{1\}$. 

\medskip 
Below we will adopt a variation of Ratcliffe's and Katok's terminology of {\em fundamental regions/domains} but impose the local finiteness condition from the beginning.  

\begin{defn}
1.  A subset $U$ of a topological space is called  an {\em open domain} (or a {\em regular open subset})  if $U$ is the interior of its 
closure. 

2. A subset $V$ of a topological space is called  a {\em closed domain} (or a {\em regular closed subset})  if $V$ is  the closure of its interior. 
\end{defn}

\begin{defn}
Suppose that $G\times X\to X$ is a proper  action of a discrete  group on a topological space $X$. 

1. An open subset $R\subset X$ is an {\em open fundamental region} for this action if the following hold:

\begin{enumerate}
\item $G\cdot \ol{R}=X$. 

\item $gR\cap R\ne \emptyset$ if and only if $g=1$. 

\item For every compact subset $K\subset X$, the transporter set $(\ol{R}|K)_G$ is finite, i.e. the family $\{g\ol{R}\}_{g\in G}$ of subsets in $X$  is {\em locally finite}. 

\end{enumerate}

2. A closed subset $F\subset X$ is a {\em closed fundamental domain} if $F$ is a closed domain in $X$ and 

\begin{enumerate}
\item $G\cdot F=X$. 

\item $g\int(F)\cap \int(F)\ne \emptyset$ if and only if $g=1$. 

\item For every compact subset $K\subset X$, the transporter set $(F|K)_G$ is finite. 
\end{enumerate}  \end{defn}

Below we describe the most common construction of open fundamental regions,  {\em open Dirichlet domains} and their variations. 

\medskip
Given an isometric metrically proper action of a discrete group $G$ on a metric space $(X,d)$  and a point $x\in X$, 
one defines the {\em open Dirichlet domain} of the action as 
$$
D_x=\{y\in X: d(y,x)< d(y, gx)\quad \forall g\in G \setminus G_x\}. 
$$   
Thus, $D_x\subset \hat{D}_x$ (see \eqref{eq:hat}). 
It is clear from the construction that for every $g\in G \setminus G_x$,
$$
g D_x \cap D_x=\emptyset 
$$
and $gD_x=D_x$ for all $g\in G_x$. In order to have a chance to get a fundamental region using open Dirichlet domains one has to assume that $G_x=\{1\}$.

\begin{rem}
Suppose that $G$ is countable, $X$ is a complete metric space, $G\times X\to X$ is a continuous action and fixed point sets in $X$ of 
nontrivial elements of $G$ are nowhere dense. Then Baire's Theorem implies existence of $x\in X$ such that $G_x=\{1\}$. For instance, if $X$ is a connected topological manifold, $G$ is discrete and acts effectively and properly on $X$.  Then the fixed-point set of each nontrivial element of $G$ has empty interior, see \cite{Newman}.   
 \end{rem}

More generally, given a closed discrete subset $E\subset X$, one defines the {\em Voronoi tessellation} 
${\mathcal V}_E$ of $X$ corresponding to $E$. The open/closed tiles of the tessellation are the subsets $V_x, \hat{V}_x, x\in E$, of $X$ defined as  
$$
V_x=\{y\in X: d(y,x)< d(y, x') \quad \forall x'\in E \setminus \{x\} \}, 
$$ 
$$
\hat V_x=\{y\in X: d(y,x)\le  d(y, x') \quad \forall x'\in E \setminus \{x\} \}.  
$$
The point $x$ is the {\em center} of the tiles $V_x, \hat V_x$. Each $\hat V_x$ is closed in $X$ (as the intersection of closed subsets).   
The {\em open tile} $V_x$ need not be an open subset of $X$ (the intersection of open subsets need not be open).  A sufficient condition is that $E\subset X$ is {\em metrically proper}, see Section \ref{sec:inmetrics}. %, i.e. every metric ball in $X$ intersects $E$ in a finite subset. 

\begin{lemma}\label{lem:loc-fin} 
If $E\subset X$ is metrically proper then each open tile $V_x$ of ${\mathcal V}_E$  is an open subset of $X$ and the collection of tiles $\hat V_x, x\in E$, is locally finite. 
\end{lemma}
\proof 1. Take $y\in V_x$. In view of metric properness of $E$, the function 
$$
d(y, \cdot)-d(y,x): E\setminus \{x\}\to \R_+ 
$$ 
attains its positive minimum $R$ at some $x' \in E\setminus \{x\}$. Then $B(y,R/2)\subset V_x$. 

2. Consider a unit ball $B=B(z,1)\subset X$. Suppose that $\hat V_x\cap B\ne \emptyset$ for some $x\in E$. Then, whenever $\hat V_y\cap B\ne \emptyset, y\in E$, $d(y,z)\le d(x,z)+1$. By the metric properness of $E$, the number of such points $y\in E$ is finite. 
\qed 

\medskip 
 The key issue that we will have to deal with is that, even if $E$ is metrically proper, the closed tile $\hat V_x$ is not necessarily the closure of the open tile $V_x$. Moreover, in general, the bisectors
$$
Bis(x,z)=\{y\in X: d(y,x)=d(y,z)\}
$$
may have nonempty interior in $X$. This happens, for instance, in the case of metric graphs. 

\begin{example}
Consider the space $X$ which is the union of two coordinate lines in $\R^2$, with the induced path-metric $d$, i.e. the restriction of the $\ell_1$-metric from $\R^2$. Thus, $(X,d)$ is a complete geodesic metric space. Let $G=\Z_2$, whose generator $g$ acts on $X$ by restriction of the antipodal map $(x,y)\mapsto (-x,-y)$ on 
$\R^2$. The group $G$ has unique fixed point in $X$, namely the origin ${\mathbf 0}=(0,0)$. For every point 
$p\in X\setminus \{{\mathbf 0}\}$ the closed Dirichlet domain $\hat D_p$ is the union of three coordinate rays, while $D_p$ consists of just one open coordinate ray. In particular, $\delta D_p=\hat D_p\setminus D_p$ is a coordinate line and, thus, is not contained in the boundary of $D_p$ (which is  the singleton $\{{\mathbf 0}\}$). The interior of  $\hat D_p$ is $\hat D_p \setminus \{{\mathbf 0}\}$, hence, 
$$
\int \hat D_p \cap g(\int \hat D_p)
$$
is nonempty and equals a coordinate line minus the origin. In particular, the interior of any closed Dirichlet domain  cannot be a fundamental region. Note also that the closure of $D_p$ is the 
closed coordinate ray containing $p$, which implies that $G \ol{D}_p\ne X$ (it misses two open coordinate rays).   
Of course, in this example one can take a suitable open subset of $\int \hat D_p$ (the union of two open rays) as a fundamental region. However, it cannot be chosen to be connected. Thus, connectedness of fundamental regions (as required by Ratcliffe's definition of a fundamental domain) is an unreasonable requirement in the setting of general complete geodesic metric spaces. 
\end{example}

Below we discuss some basic properties of Voronoi tiles. 

\begin{lem}
Let $\phi: X\to (0,\infty)$ be an $L$-Lipschitz function for some $L\le 1/2$ and let $E\subset X$ be such that for every $x\in X$, the open ball $B(x, \phi(x))$ has nonempty intersection with $E$. Then for every $x\in X$, $\hat V_x\subset  B(x, 2 \phi(x))$. 
\end{lem}
\proof Take $y\in \hat V_x$. Then there exists $z\in E$ such that $d(y,z)< \phi(y)$. Since  $y\in \hat V_x$, $d(x,y)\le d(z,y)<\phi(y)$. 
By the $L$-Lipschitz property of $\phi$, we have $\phi(y)\le \phi(x) +L d(x,y)$, implying
$$
(1-L)d(x,y) <  \phi(x)
$$
and $d(x,y)< \frac{1}{1-L} \phi(x) \le 2\phi(x)$. Therefore, $y\in  B(x, 2 \phi(x))$. \qed 

\medskip
Suppose that $(X,d)$ is a metric space, $G\times X\to X$ an isometric properly discontinuous action, $E\subset X$ is a metrically proper $G$-invariant subset, $S=q(E)\subset X/G$ is the image of $E$ under the quotient map $q: X\to X/G$. We then have two Voronoi tessellations $\mathcal V_E$ (of $X$) and $\mathcal V_S$ (of $X/G$ equipped with the metric $d_G$). 

\begin{lemma}\label{lem:onto}
For every closed and every open Voronoi tile $\hat V_x, V_x, x\in E$, we have $q(\hat V_x)=\hat V_{[x]}$ and 
$q(V_x)=V_{[x]}$. 
\end{lemma}
\proof The statement is a direct consequence of definitions of Voronoi tiles and the metric $d_G$. \qed 

\medskip 
Our next goal is to find a condition on $E$ that ensures injectivity of the restriction of $q$ to each $\hat V_x$. Recall that in Section \ref{sec:inmetrics}, given an isometric properly discontinuous action $G\times X\to X$, we defined a function $\rho: X/G\to \R_+$ (as well as $\rho: X\to \R_+$). 

\begin{lem}\label{lem:injective}
Suppose now that $G$ is a group acting freely, isometrically and properly discontinuously on $(X,d)$. 
Suppose that $E\subset X$ is a $G$-invariant closed discrete subset such that $B(x, \frac{1}{4}\rho(x))\cap E\ne \emptyset$ for 
all $x\in X$. Then the quotient map $q: X\to X/G$ is injective on each closed tile $\hat V_x, x\in E$. In other words, 
$$
\hat V_x\cap \hat V_{gx}=\emptyset, \ \forall x\in E, \ \forall g \in G\setminus \{1\}. 
$$
\end{lem}
\proof For every $y\in \hat V_x$ there exists $x'\in E$ at distance $< r=\frac{1}{4}\rho(y)$ from $y$. Since  
$y\in \hat V_x$, we have $d(x,y)< r$ and, therefore, $d(x,x')< 2r$. 
Using the inequality  
$$
|\rho(x) - \rho(y)|< 2r,
$$
we get $d(x, gx)\ge \rho(x)>  \rho(y) - 2r =2r$ for all $g \in G\setminus \{1\}$. By the triangle inequality, $d(x,y)< r$ implies that $d(y, gx)> r$. Thus, $y\notin \hat V_{gx}$. \qed 

\medskip
We will construct subsets $E\subset X$ satisfying the assumptions of Lemma \ref{lem:injective} in Lemma \ref{lem:max-net}. 

%We will construct a certain subset $S\subset X/G$ such that $B([x], \frac{1}{4}\rho(x))\cap S\ne \emptyset$ for all $[x]\in X/G$. 

\medskip 
A subset $A$ of a geodesic metric space $(X,d)$ is {\em starlike} with respect to a point $a\in A$ if for each $x\in A$ every geodesic segment $ax$ is contained in $A$.

\begin{lemma}\label{lem:starlike}
Suppose that $(X,d)$ is a geodesic space and $\mathcal V_E$ be the Voronoi tessellation corresponding to a metrically proper subset $E\subset X$.  Then each tile $V_x, \hat V_x$ of $\mathcal V_E$ 
 is starlike with respect to its center. 
\end{lemma}
\proof The proof is essentially the same as the one in \cite[Theorem 6.6.13]{Ratcliffe}. Take a point $z\in \hat{V}_x$ and let 
$c: [0,T]\to X$ be a geodesic connecting $x$ to $z$. Then for each $t\in [0,T]$ and $y\in E\setminus \{x\}$, we get (by the triangle inequality) 
$$
d(x, c(t))=t=T - d(c(t), z)=d(x,z)- d(c(t), z)\le d(y,z)- d(c(t), z) \le d(y, c(t)).$$
Hence $c(t)\in \hat V_x$ and, therefore, $\hat V_x$ is starlike with respect to $x$. The same argument works for $V_x$. \qed 
 
\medskip 
The basic examples of Voronoi tessellations are when $(X,d)$ is a Euclidean or a  real-hyperbolic space; in these cases  Voronoi tiles 
(and, hence, their intersections) are convex. This need not be the case in general even when one works with, say, {\em complex-hyperbolic spaces} (see e.g. \cite{Goldman}). Below we will see that some kind of convexity still holds in the case of Voronoi tessellations of Gromov-hyperbolic spaces. 

Recall that a subset $Y$ of a geodesic metric space $(X,d)$ is called {\em $\la$-quasiconvex} if every geodesic segment $xy$ with the end-points in $Y$ is contained in the closed $\lambda$-neighborhood of $Y$, i.e. $d(z, Y)\le \la$ for all $z\in xy$. 

\begin{cor}
Suppose, additionally, that $(X,d)$ is $\delta$-hyperbolic. Then:

1. For every Voronoi tessellation $\mathcal V_E$ of $X$, each   tile $V_x, \hat V_x$ is $\delta$-quasiconvex. 

2. Each bisector $Bis(x,y)$ in $X$ is  
%Intersections of pairs of tiles $\hat V_{x}\cap \hat V_{y}$ are 
$2\delta$-quasiconvex. 
\end{cor}
\proof 1. This is a direct consequence of Lemma \ref{lem:starlike} and the definition of $\delta$-hyperbolicity via slimness of geodesic triangles. 

2. Take $E=\{x,y\}$ and the corresponding Voronoi tessellation of $X$ with just two closed tiles, $\hat V_{x}, \hat V_{y}$.  
Then $Bis(x,y)= \hat V_{x}\cap  \hat V_{y}$. Suppose that points $p, q$ belong to  $Bis(x,y)$. 
Consider a point $z$ on a geodesic $pq$ in $X$. 
By Part 1, there exist points $x'\in \hat V_x,  y'\in \hat V_y$ within distance $\delta$ from $z$. In particular, $d(x',y')\le 2\delta$. Since the geodesic $x'y'$ connects $\hat V_x, \hat V_y$, by continuity of the function $d(x, \cdot )- d(y, \cdot)$, there exists $z'\in x'y'\in Bis(x,y)$. 
By the triangle inequality, $d(z, z')\le 2\delta$.  \qed 

\medskip
Note that the proof of Lemma \ref{lem:starlike} also shows that for each $z\in \hat V_x$ and $t< T$, either we get the strict inequality  
$d(x, c(t))< d(y, c(t))$, or $c(t)$ belongs to a geodesic $yz$. If the former case occurs for all $y\in E\setminus \{x\}$, we conclude that 
we get $c(t)\in V_x$. In particular, in that case, $\hat V_x$ is the closure of $V_x$. In order to rule out the second possibility ($c(t)$ belongs to a geodesic $yz$) one has to impose extra restrictions. 

\begin{defn}
A geodesic metric space $(X,d)$ has {\em nonbranching geodesics} if each geodesic $c: I\to (X,d)$ is uniquely determined by its restriction to a nonempty open subset of the interval $I$.  
\end{defn}

For instance, geodesics in Riemannian manifolds and, more generally, manifolds with smooth Finsler metrics and Alexandrov spaces  satisfy this property. 

\begin{cor}\label{cor:nonbranching1}
Suppose that $(X,d)$ is a metric space with nonbranching geodesics. Then for each Voronoi tessellation $\mathcal V_E$ of $X$ and every $x\in E, z\in \hat{V}_x$ and geodesic $c: [0,T]\to xz$ connecting $x$ to $z$, one has $c(t)\in V_x$, $t< T$. In particular, $\hat V_x$ is the closure of $V_x$. Moreover, $V_x$ is an open domain in $X$. 
\end{cor}
\proof The proof is the same as the one in \cite[Theorem 6.6.13]{Ratcliffe}. Suppose that $c(t)\notin V_x$ for certain $t<T$. Then 
for all $s\in [t,T]$ we have $c(s)\in \delta V_x= \hat{V}_x \setminus V_x$. Due to the local finiteness of $\mathcal V_E$, there exist  
$y\in E \setminus \{x\}$ and $t'$, $t< t' <T$, such that $c(s)\in \hat V_y$ for all $s\in [t', T]\subset [t,T]$. Therefore, for all $s\in [t', T]$ we get 
$$
s=d(x,c(s))= d(y, c(s)),
$$
and, thus, $d(y, c(s))+ d(c(s),z)= d(y,z)$. In other words, the concatenation of geodesics $yc(t')\star c(t')z$ is a geodesic $\gamma$ in 
$(X,d)$. For the second  segment $c(t')z$ of this geodesic we will take the restriction of $c$ to $[t', T]$. Since $y\ne x$, geodesics $c$ and $\gamma$ have distinct images; on the other hand, they agree on the open subinterval $(t',T)$. This contradicts the nonbranching assumption. The proof that $V_x$ is an open domain in $X$ is similar and we omit it. \qed   

\begin{cor}
[See Theorem 6.6.13 in \cite{Ratcliffe}] 
Suppose that $(X,d)$ is a metric space with nonbranching geodesics, $G\times X\to X$ is an isometric metrically proper action. Then for every $x\in X$ with $G_x=\{1\}$ the open Dirichlet domain $D_x$ is a connected open fundamental region for the $G$-action on $X$. Moreover, $\hat D_x$ is the closure of $D_x$ and the interior of $\hat D_x$ is precisely $D_x$. 
\end{cor}
\proof We apply Corollary \ref{cor:nonbranching1} to the Voronoi tessellation $\mathcal V_{Gx}$. Corollary \ref{cor:nonbranching1} implies that $\hat{D}_x$ is the closure of $D_x$. Connectedness of $D_x$ is clear from the same corollary. The fact that $G \hat D_x=X$ follows from Proposition \ref{prop:Dirichlet domain}. It remains to prove that  the interior of $\hat D_x$ is $D_x$. Take $z\in \hat D_x\cap \hat D_y$, where $y\in Gx \setminus  \{x\}$. By applying Corollary \ref{cor:nonbranching1} to the Voronoi tile $\hat D_y$ we see that 
$z$ belongs to the closure of $D_y$. Since the latter is disjoint from $\hat D_x$, we conclude that $z\notin \int \hat D_x$.  \qed

\begin{question}
Suppose that $M$ is a connected topological manifold. Does $M$ admit a complete geodesic metric with nonbranching geodesics? 
\end{question}

In the rest of the section we will prove existence of connected closed fundamental domains for free properly discontinuous actions on geodesic metric spaces by using Voronoi tessellations more general than the ones given by the Dirichlet construction. (More precisely, we will use 
$\mathcal V_{E}$ for some $G$-invariant closed discrete subset $E$ of $X$.) 
%We do not know how to find connected open fundamental regions under this assumption. 
In what follows, $(X,d)$ is a separable geodesic complete  metric space. 

\begin{lemma}
Fix a point $y\in X$. Then there is a $G_\delta$-subset $X_y\subset X$ consisting of points $x$ such that $Bis(x,y)$ has empty interior. 
\end{lemma}
\proof First of all, we prove that for every $z\in X$ the subset $\{x\in X: z\notin Bis(x,y)\}$ is open and dense in $X$. Openness is clear. 
To prove denseness, take a sequence of points $x_n\in xz\setminus \{x\}$ converging to $x$ and note that 
for every $x_n$, $z\notin Bis(x_n,y)$.  

Now, we take a dense countable subset $Z\subset X$ (here we use the separability assumption). Define
$$
X_y=\{x: \forall z\in Z, z\notin Bis(x,y)\}. 
$$
Then $X_y$ is the intersection of a countable family of open and dense subsets, i.e. is a $G_\delta$-subset of $X$. Let us prove that for 
every $x\in X_y$ the bisector $Bis(x,y)$ has empty interior. Take $w\in Bis(x,y)$ and a sequence $z_n\in Z$ converging to $w$. 
Then, by the definition of the set $X_y$, for all $n$, $z_n\notin Bis(x,y)$. Hence, $w$ cannot belong to the interior of $Bis(x,y)$. 
\qed 

Since the intersection of a countable family of $G_\delta$-subsets is again a $G_\delta$-subset, we obtain:

\begin{cor}\label{cor:empty-interior} 
Let $Y\subset X$ be a countable subset. Define $X_Y=\bigcap_{y\in Y} X_y$. Then $X_Y$ is a $G_\delta$-subset of $X$. For every $x\in X_Y$ and $y\in Y$ the bisector $Bis(x,y)$ has empty interior. 
\end{cor}

\begin{lemma}\label{lem:max-net} 
Suppose that $(M, d)$ is a separable metric space, $\phi: M\to (0,\infty)$ a continuous function. Then there exists a countable discrete and closed subset $C\subset M$ such that for all $z\in M$, $B(z, \phi(z))\cap C$ is nonempty. 
\end{lemma}
\proof Take a maximal subset $S$ of $M$ satisfying the property that for all distinct $x, y\in S$, 
$$
d(x,y)\ge \frac{1}{2} \min( \phi(x), \phi(y)). 
$$ 
Separability of $M$ ensures that such a subset is countable. Suppose that $(x_n)$ is a sequence of distinct elements in $S$ converging to $x\in M$. By continuity of $\phi$, we have that 
$$
\phi(x_n) > \eps=\frac{1}{2} \phi(x)>0
$$ 
for all sufficiently large $n$. %Without loss of generality, we may assume that all the points $x_n$ belong to the ball $B(x, \eps)$. 
We have (for all $m\ne n$) 
$$
d(x_m, x_n)\ge  \frac{1}{2} \min( \phi(x), \phi(y))>  \eps.
$$
But then the sequence $(x_n)$ cannot converge. 

Lastly, take $x\in M$. Suppose that $B(x, \phi(x))\cap S=\emptyset$. Then $d(x, y)\ge \phi(x)$ for all $y\in S$, which implies  
$$
d(x,y)\ge  \min( \phi(x), \phi(y)). 
$$
Hence, $S\cup \{x\}$ still satisfies the inequality defining $S$. This contradicts the maximality of $S$. \qed 

\begin{addendum}
Assume, additionally, that $(M,d)$ is a proper metric space.  
Then the collection of Voronoi tiles $\hat V_x, x\in C$, associated with the subset $C$ in the lemma, is locally finite.  
\end{addendum}
\proof This follows from the fact that each metric ball $B(x,R)$ in $M$ contains only finitely many points from $C$, see Lemma 
\ref{lem:loc-fin}. \qed 

We next perturb the subset $C$ in Lemma \ref{lem:max-net} to a new subset $C'$ which satisfies essentially the same properties as $C$ but also has empty interior of $\delta V_{x'}:=\hat V_{x'}\setminus V_{x'}$ for every $x'\in C'$:  

 \begin{lemma}\label{lem:perturb}
 Suppose that $(M,d)$ is a complete separable geodesic metric space. 
 Then there exists a countable discrete and closed subset $C'\subset M$ such that for all $z\in M$, $B(z, 2\phi(z))\cap C$ is nonempty 
 and  for each $x'\in C'$, $\delta V_{x'}$ has empty interior. 
 \end{lemma}
\proof We will take the subset $C$ constructed in the previous lemma and perturb it inductively using Corollary \ref{cor:empty-interior} so that 
the bisectors $Bis(x',y')\subset M$ (for distinct points $x', y'\in C'$) have empty interior. Since  the complements   
$\delta V_{x'}$ are contained in the union of bisectors $Bis(x',y'), y'\in C'\setminus \{x'\}$, it will follows that each $\delta V_{x'}$ has empty interior. 

In order to construct the perturbation, using continuity and nonvanishing of $\phi$, 
for every $x\in C$ we find $\eps_x>$ such that:

1. For all $z\in M$ satisfying $x\in B(z,\phi(z))$ we have $\eps_x< \phi(z)$. 

2. $\lim_{n\to\infty} \eps_{x_n} =0$ for some enumeration of the set $C$.  

3. The metric balls $B(x, \eps_x)$, $x\in C$, are pairwise disjoint. 

Now, we replace each $x\in C$ with some $x'\in B(x, \eps_x)$; the resulting 
subset $C'$ of $M$ is still countable. Furthermore, for every $z\in M$, there exists $x'\in C'$ as above such that 
$d(x',z)< 2\phi(z)$. If $C'$ fails to be closed and discrete, there exists a sequence $x'_n$ of distinct elements of $C'$ converging to 
some $x\in M$.  By (2), 
$$
\lim_{n\to\infty} d(x_n, x'_n)=0, \quad x_n\in C.  
$$
Then the sequence $(x_n)$ also converges to $x$ and, by (3), all the points $x_n$ are distinct. This is a contradiction. \qed 

\begin{lemma}\label{lem:closures-cover} 
Suppose, additionally, that $(M,d)$ is a proper metric space (separability is automatic in this case). 
Then for every point $y\in M$ there exists $x\in C'$ such that $y\in \ol{V}_x$. 
%belongs to the closure of the interior 
\end{lemma}
\proof  By the local finiteness of the Voronoi tessellation, for each point $y\in M$ there is a neighborhood $U$ of $y$ which intersects only finitely many closed tiles $\hat V_{x_i}, i=1,...,n$. Suppose that 
$$
y\in \delta V_{x_i}, i=1,...,n. 
$$
Then, by shrinking the neighborhood $U$ further, we can assume that $U\subset \delta V_{x_1}\cup ... \cup \delta V_{x_n}$. But this contradicts the fact that each $ \delta V_{x_i}$ has empty interior. \qed

\medskip
We now consider the situation when $(X,d)$ is a proper geodesic metric space, $G\times X\to X$ is a free, isometric, properly discontinuous action. We form the quotient space $(M,d_M)=(X/G, d_G)$. This space is again proper and geodesic, see Lemma \ref{lem:quot-space}. Take the function $\phi=\frac{1}{16}\rho: M\to \R_+$, where $\rho$ is defined via the $G$-action on $X$ as in   
Section \ref{sec:inmetrics}. Using Lemma \ref{lem:perturb} we find a suitable (countable) metrically proper subset $C'\subset M$. 
Let $E\subset X$ denote $q^{-1}(C')$, where 
$q: X\to X/G=M$ is the quotient map. We obtain two Voronoi tessellations, $\mathcal V_E, \mathcal V_{C'}$ of $X$ and $M$ respectively. 
According to Lemmata \ref{lem:onto} and \ref{lem:injective}, 
the map $q$ sends each closed tile $\hat V_x, x\in E$, of $\mathcal V_E$ homeomorphically to the closed tile 
$\hat V_{[x]}$ of $\mathcal V_{C'}$. By applying Lemma \ref{lem:closures-cover}, we obtain:

\begin{cor}\label{cor:closures-cover} 
For every point $y\in X$ there exists $x\in E$ such that $y\in \ol{V}_x$. 
\end{cor}  

%From each orbit $Gx, x\in E$, we pick one point and, thus, form a subset $S\subset E$ consisting of such points. 
We are now ready to prove the main theorem of this section: 

\begin{thm}\label{thm:region} 
Suppose that $(X,d)$ is a proper geodesic metric space, $G\times X\to X$ is a free, isometric, properly discontinuous action. Then 
this action admits an open fundamental region $R$ and a closed connected fundamental domain $F$. 
\end{thm}
\proof We define a $G$-invariant subset $E\subset X$ as above and pick its subset $S\subset E$ intersecting each $G$-orbit in $E$ in exactly one point. Take
$$
R=\bigcup_{x\in S} V_x. 
$$
We claim that $R$ is an open fundamental region for the $G$-action on $X$. First of all, since $G$ preserves the open and closed tiles of the Voronoi tessellation $\mathcal V_E$, we have that for every $x\in S$, $gV_x\cap V_x\ne \emptyset$ if and only if $gx=x$, i.e. $g=1$ (since $g$ acts freely on $X$). Next, the tiling  $\mathcal V_E$ is locally finite (since the subset $E$ is properly embedded in $X$). Hence, the collection of closures $\ol{V}_x, x\in E$, is locally finite as well. Lastly, by Corollary \ref{cor:closures-cover}, each point 
$y\in X$ belongs to some  $\ol{V}_x, x\in E$. Then there exists $z\in S$ such that $g(z)=x$ and, hence, $g \ol{V}_z= \ol{V}_x$. 
Thus, all properties of an open fundamental region are satisfied by $R$. 

\medskip
We construct a closed connected fundamental domain 
$F$ by modifying the construction of $R$ above. The main issue is that for a random choice of $S$, the region $R$ need not have 
connected closure. We will use Lemma \ref{lem:tree} to choose $S$ more carefully. We let $\Gamma$ denote the {\em incidence graph} of 
the closed cover $\{\ol{V}_x: x\in E\}$ of $X$: Vertices of $\Gamma$ are the elements of $E$ and we connect distinct points $x, y\in E$ by an edge if and only if $\ol{V}_x\cap \ol{V}_y\ne \emptyset$. Since $X$ is connected, the graph $\Gamma$ is connected as well. By the construction, the graph $\Ga$ is simplicial and the group   $G$ acts on $\Ga$ by simplicial automorphisms. We let $\Ga'$ denote the barycentric subdivision of $\Ga$ and $\Phi\subset \Ga'$ a subtree given by  Lemma \ref{lem:tree}. We let $S\subset E$ denote the intersection of $E$ with the vertex set of $\Phi$ (i.e. the subset of vertices of $\Phi$ which are vertices of $\Ga$). According to 
 Lemma \ref{lem:tree},  each $G$-orbit in $E$ intersects $S$ in exactly one point. Thus, the open subset $R\subset X$ defined as above using $S$ is an open fundamental region for the $G$-action on $X$.  We let 
$$
F:= \bigcup_{x\in S} \ol{V}_x. 
$$ 
Let us verify connectedness of $F$. First of all, each $\ol{V}_x$ is connected since $V_x$ is connected. For any two points $x, y\in S$ 
there exists a vertex-path $x_1 x_2 ... x_n$ in $\Ga$ connecting these vertices ($x_1=x, x_n=y$) such that each vertex of this path is in $S$ (this follows from the fact that the graph $\Phi$ is connected). The union
$$
\ol{V}_{x_1}\cup... \cup \ol{V}_{x_n}
$$
is connected since each intersection $\ol{V}_{x_i}\cap \ol{V}_{x_{i+1}}$ is nonempty. Thus, $F$ is connected. The fact that $F$ is closed follows from the fact that each $\ol{V}_x$ is closed and the family $\{\ol{V}_{x}: x\in S\}$ is locally finite in $X$. The subset 
$$
R= \bigcup_{x\in S} {V}_x
$$
is open in $X$ and dense in $F$. Hence, $F$ is a closed domain in $X$. \qed

\medskip 
As an application we will prove existence of open fundamental regions and closed connected fundamental domains for free properly discontinuous group actions on a certain class of topological spaces, cf. \cite{MO3}.

\begin{theorem}
Suppose that $X$ is a 2nd countable, connected and locally connected locally compact Hausdorff topological space. 
Suppose that $G\times X\to X$ is a free proper action of a discrete countable group. Then this action admits an open fundamental region and a closed connected fundamental domain. 
\end{theorem}  
\proof In Theorem \ref{thm:geodesic-metrization} we constructed a complete $G$-invariant geodesic metric $d$ on $X$. The metric space $(X,d)$ is necessarily proper, since $X$ is locally compact. Using Theorem \ref{thm:region} we find an open fundamental region $R$ and a closed connected fundamental domain $F$ for the $G$-action on $X$. 
\qed 
  
 % \newpage 

\end{document}